\documentclass[12pt]{article}

\usepackage{amsfonts, amsmath, amsthm, amssymb, tikz, url}
\usepackage[in]{fullpage}
\usepackage[english]{babel}
\usepackage[utf8]{inputenc}
\usepackage{authblk}
\usepackage{float}
\usepackage{subfigure}

\usetikzlibrary{shapes.geometric}

\newtheorem{theorem}{Theorem}
\newtheorem{claim}{Claim}

\title{Specialization, generalization, and a new proof of Viviani's Theorem}
\author{Addie Armstrong}
\author{Dan McQuillan}
\affil{Department of Mathematics, Norwich University, Northfield VT 05663}

\begin{document}
\maketitle

\begin{quote}
\textit{The first step deals with a particular case which is not only especially accessible but also especially useful; we can appropriately call it a} leading particular case\textit{: it leads the way to a general solution.}
\end{quote}
\vspace{-0.7cm}
\begin{center}
George Polya
\end{center}

One of the primary goals of mathematics instruction is to communicate a sense of how mathematics is done. In his book \cite{Polya}, George Polya illustrates the idea behind his quote with a well known proof of the Central Angle Theorem; the special case where a diameter is used to form one of the angles being easy to prove and (more interestingly) quickly leading to a complete proof for the general case. Thus, a carefully chosen special case can truly be seen as being at the heart of a full proof, which is completed by what Polya calls {\it superposition of particular cases}. 

The goal of this note is to provide an apparently new example of Polya's quoted suggestion of using a {\it leading particular case} en route to a general solution. Our proof of Viviani's theorem also highlights the value of generalization in mathematics; the result being ostensibly about an equilateral triangle, with the proof crucially using statements about {\it all} equilateral triangles. This proof is so simple that it can be used at all levels of mathematics instruction to illustrate the relationship between special examples and general statements.

\begin{theorem}[Viviani, circa 1659]
For any point $P$ on the interior or boundary of an equilateral triangle $\Delta ABC$, the sum of the distances from $P$ to each side is equal to the altitude of triangle $\Delta ABC$.
\end{theorem}

The idea of the proof is to focus on the the special case where $P$ lies on one of the sides of triangle $\Delta ABC$. Once this is proved, it is not only true for triangle $\Delta ABC$ itself, but for all equilateral triangles. Moving from the special case to the general case involves using the special case, but on a different triangle (triangle $\Delta MBN$ in figure 3). {\it Interestingly, this special case is itself proved via its own special case}---a leading particular case within a leading particular case!

\begin{proof}
We begin with a case that is completely clear, since $\Delta ABC$ is an equilateral triangle:
\begin{claim}
If $P$ lies on a vertex of $\Delta ABC$ (see Figure \ref{Case 1}) then the conclusion of Viviani's Theorem is true.
\end{claim}

\begin{figure}[H]
\begin{center}
\begin{tikzpicture}
\node[draw, thick, black ,minimum size=2.4cm,regular polygon,regular polygon sides=3] (a) {};

\filldraw[black] (a.corner 1) circle[radius=2.5pt];
\draw (a.corner 1) node[above]{$P$};
\draw[thick, blue, dashed] (a.corner 1) --(0,-0.6);
\draw (a.corner 2) node[left] {$A$};
\draw (a.corner 3) node[right]{$C$};
\end{tikzpicture}
\end{center}
\caption{$P$ is a vertex of $\Delta ABC$}
\label{Case 1}
\end{figure}
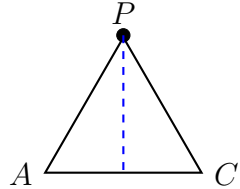

\begin{claim} If $P$ lies on an edge of $\Delta ABC$ then the conclusion of Viviani's Theorem is true.
\end{claim}

\textit{Proof of Claim 2:} Without loss of generality we may assume $P$ is on the edge $BC$, as shown in Figure \ref{Case 2}. Let $Q$ be the point of intersection of $AB$ and the line through $P$ parallel to $AC$.
Since $\Delta QBP$ is an equilateral triangle, we can use the result of Claim 1 on triangle $\Delta QBP$.
Since $QP$ is parallel to $AC$, the distance from $P$ to $AC$ is the same as the distance from any other point on $QP$ to $AC$. This proves Claim 2. 

\begin{figure}[H]
\begin{center}

\begin{subfigure}[]{}
\begin{tikzpicture}
\node[draw, thick, black ,minimum size=2.4cm,regular polygon,regular polygon sides=3] (-1,0) (b) {};

\draw (0.66,0) node[right]{$P$};

\draw [thick, black] (0.66,0)--(-0.66,0);
\draw (-0.66,0) node[left]{$Q$};

\filldraw[black] (0.66,0) circle[radius=2.5pt];

\draw (a.corner 1) node[above]{$B$};
\draw[thick, black, dotted] (0,0) --(0,-0.6);
\draw[thick, black, dashed] (a.corner 1)--(0,0);
\draw (a.corner 2) node[left] {$A$};
\draw (a.corner 3) node[right]{$C$};
\end{tikzpicture}
\end{subfigure} \hspace{0.5cm}
\begin{subfigure}[]{}
\begin{tikzpicture}

\node[draw, thick, black ,minimum size=2.4cm,regular polygon,regular polygon sides=3] (a) {};

\draw (0.66, 0) node[right] {$P$};

\draw[thick, blue, dashed] (0.66,0)--(-0.33,0.6);
\draw[thick, blue, dotted] (0.66,0)--(0.66,-0.62);

\draw [thick, black] (0.66,0)--(-0.66,0);
\draw (-0.66,0) node[left]{$Q$};

\filldraw[black] (0.66,0) circle[radius=2.5pt];

\draw (a.corner 1) node[above]{$B$};
\draw[thick, black, dotted] (0,0) --(0,-0.6);
\draw[thick, black, dashed] (a.corner 1)--(0,0);
\draw (a.corner 2) node[left] {$A$};
\draw (a.corner 3) node[right]{$C$};
\end{tikzpicture}
\end{subfigure}
\end{center}
\caption{$P$ is on an edge of $\Delta ABC$}
\label{Case 2}

\end{figure}
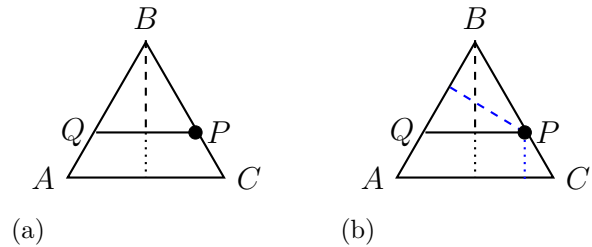

\begin{claim}
If $P$ lies in the interior of $\Delta ABC$ then the conclusion of Viviani's Theorem is true.
\end{claim}

\textit{Proof of Claim 3:} Let $M$ and $N$ be points on $AB$ and $BC$, respectively, such that MN is parallel to AC and $P$ lies on $MN$. We use Claim 2 on the triangle $\Delta MBN$. The result now follows from the fact that the distance from $P$ to $AC$ is equal to the distance from any point on $MN$ to $AC$. This proves the claim and the theorem.

\end{proof}

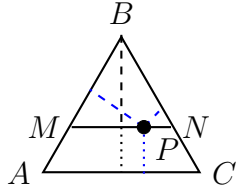
\begin{figure}[H]
\begin{center}
\begin{tikzpicture}
\node[draw, thick, black ,minimum size=2.4cm,regular polygon,regular polygon sides=3] (a) {};

\draw (0.3, 0) node[below right] {$P$};

\draw[thick, blue, dashed] (.3,0)--(-0.42,0.5);
\draw[thick, blue, dotted] (.3,0)--(0.3,-0.62); 
\draw[thick, blue, dashed] (0.3,0)--(0.5,0.21);

\draw [thick, black] (.66,0)--(-0.66,0);
\draw (-0.66,0) node[left]{$M$};
\draw (0.66,0) node[right]{$N$};

\filldraw[black] (0.3,0) circle[radius=2.5pt];

\draw (a.corner 1) node[above]{$B$};
\draw[thick, black, dotted] (0,0) --(0,-0.6);
\draw[thick, black, dashed] (a.corner 1)--(0,0);
\draw (a.corner 2) node[left] {$A$};
\draw (a.corner 3) node[right]{$C$};
\end{tikzpicture}
\end{center}
\caption{$P$ is in the interior of $\Delta ABC$}
\label{Case 3}
\end{figure}

While many beautiful proofs of Viviani's Theorem are available (see for example \cite{Kawasaki} \cite{Tanton} \cite{Wolf}), our proof has the special feature that it highlights the relationship between specialization and generalization. This challenges one's imagination to move from triangle to triangle, illustrating the central role that generalization has in mathematics. 

Furthermore, related theorems can be easily proved using this framework. In his video on Vivani's Theorem, \cite{Tanton}, James Tanton discusses similar results for points outside the triangle.
This leads to the next claim.

\begin{claim}
Extend the sides of triangle $\Delta ABC$ to full lines. Let $P$ be in one of the infinite regions bounded by an edge of $\Delta ABC$. Let $x$ be the distance between $P$ and the extension of the bounding edge. Let $y$ and $z$ be the distances from $P$ to the extensions of the other edges. Then $y+z-x$ is equal to the altitude of triangle $\Delta ABC$.
\end{claim}

\begin{figure}[H]
\begin{center}
\begin{tikzpicture}
\node[draw, thick, black ,minimum size=2.4cm,regular polygon,regular polygon sides=3] (a) {};

\draw (1.7,-2.5) node[below left] {$P$};

\draw[thick, blue, dashed] (1.7,-2.5)--(-1.2,-0.9);
\draw[thick, blue, dotted] (1.7,-.6)--(1.7,-2.5);
\draw (1.7, -1.2) node[right] {$x$};
\draw[thick, blue, dashed] (1.7,-2.5)--(2,-2.2);

\draw [thick, black] (2.5,-2.5)--(-2.5,-2.5);

\filldraw[black] (1.7,-2.5) circle[radius=2.5pt];

\draw (a.corner 1) node[left]{$B$};
\draw[thick, black, dotted] (0,-.6) --(0,-2.5);
\draw[thick, black, dashed] (a.corner 1)--(0,-.6);
\draw (a.corner 2) node[above left] {$A$};
\draw (a.corner 3) node[above right]{$C$};

\draw[thick, black] (0.3, 1.7)--(-2.52,-3.1);
\draw[thick, black] (-.3, 1.7)--(2.52,-3.1);
\draw[thick, black] (-3,-0.6)--(3,-0.6);

\draw (-2.2, -2.6) node[ above left] {$M$};
\draw (2.2, -2.6) node[above right] {$N$};
\end{tikzpicture}
\end{center}
\caption{$P$ is in an infinite region bounded by an edge of $\Delta ABC$}
\label{Case 4}
\end{figure}
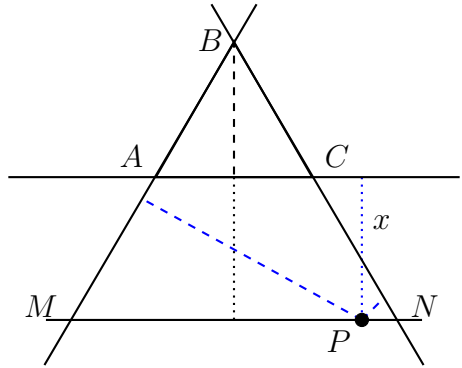

The proof of Claim 4 is identical to the proof of Claim 3 with the understanding that $AB$, $BC$, and $AC$ refer to the extensions of the sides rather than just the sides themselves. (see Figure \ref{Case 4}).

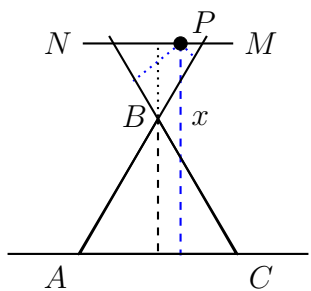
\begin{figure}[H]
\begin{center}
\begin{tikzpicture}
\node[draw, thick, black ,minimum size=2.4cm,regular polygon,regular polygon sides=3] (a) {};

\draw (0.3, 2.2) node[above right] {$P$};

\draw[thick, blue, dotted] (.3,2.2)--(-0.33,1.7);
\draw[thick, blue, dashed] (.3,2.2)--(0.3,-0.62);
\draw (0.3, 1.2) node[right] {$x$};
\draw[thick, blue, dotted] (0.3,2.2)--(0.5,2);

\draw [thick, black] (1,2.2)--(-1,2.2);
\draw (-1,2.2) node[left]{$N$};
\draw (1,2.2) node[right]{$M$};

\draw[thick, black] (0.65, 2.3)--(a.corner 2);
\draw[thick, black] (-.65, 2.3)--(a.corner 3);
\draw[thick, black] (-2,-0.6)--(2,-0.6);

\filldraw[black] (0.3,2.2) circle[radius=2.5pt];

\draw (a.corner 1) node[left]{$B$};
\draw[thick, black, dotted] (a.corner 1) --(0,2.2);
\draw[thick, black, dashed] (a.corner 1)--(0,-.6);
\draw (a.corner 2) node[below left] {$A$};
\draw (a.corner 3) node[below right]{$C$};
\end{tikzpicture}
\end{center}
\caption{$P$ is in an infinite region not bounded by an edge of $\Delta ABC$}
\label{Case 5}
\end{figure}

Finally, if $P$ is in an infinite region not bounded by an edge of triangle $\Delta ABC$, there is still a vertex, $B$ say, on the boundary of the region containing $P$ (see Figure \ref{Case 5}). Let $x$ be the distance from $P$ to the extension of $AC$, and $y$ and $z$ be defined as in Claim 4. Let $M$ and $N$ be defined as in the proof of Claim 3. A similar argument shows that $x - y -z$ is the altitude of triangle $\Delta ABC$.

For a three dimensional analog of Viviani's theorem, consider a point $P$ inside of a regular tetrahedron; from the perspective of the point $P$, imagine the floor rising until it meets $P$. Then the walls close in on $P$ until $P$ is a vertex of a smaller regular tetrahedron. Can this vague idea be turned into a proof? Use the idea of specialization to find out! Many other generalizations and variations of Viviani's theorem are possible and we encourage the interested reader to explore them.

\end{document}